# New Approach to Identify Common Eigenvalues of real matrices using Gerschgorin Theorem and Bisection method

T. D. Roopamala

Dept. of Computer science and Engg;

SJCE, Mysore, India.

S. K. Katti

Dept. of Computer science and Engg,

SJCE, Mysore, India.

*Abstract*—In this paper, a new approach is presented to determine common eigenvalues of two matrices. It is based on Gerschgorin theorem and Bisection method. The proposed approach is simple and can be useful in image processing and noise estimation.

## INTRODUCTION

Eigenvalues play vary important role in engineering applications. A vast amount of literature is also available for computing eigenvalues of a given matrix. Moreover, various numerical techniques such as bisection method, Newton Raphson method, Regula falsi method etc., are available for computing eigenvalues [3]. These methods are applied in various engineering applications. In practice, for some applications, common eigenvalues of the matrices are required. These eigenvalues can be calculated using above methods. In [5], an algorithm is presented to identify common eigenvalues of a two matrices without computing actual eigenvalues of the matrices. But, this method requires Hessenberg transformation of matrices. While going various literature survey, it is observed that except algorithm as proposed [5], no algorithm is available which can be used to identify common eigenvalues of matrices.

In practice identification of common eigenvalues are required in various image processing, control systems, and noise estimation applications. Therefore, in this paper attempt has been made to identify common eigenvalues using Gerschgorin theorem and Bisection method. The proposed approach is improvement over the bisection method for computing common eigenvalues.

In this paper, Gerschgorin circles have been drawn for two matrices. Then by selecting intersection area of two matrices, bound under which all real common eigenvalues lying are determined. This improved bound is considered as initial approximation for computing eigenvalues of two system matrices using Bisection method. These are compared with the Bisection method with approximate (trial) approximation.

## II. GERSCHGORIN THEOREM [1]

For a given matrix *A* of order ($n \times n$), let $P_k$ be the sum of the moduli of the elements along the $k^{th}$ row excluding the diagonal elements $a_{kk}$. Then every eigenvalues of *A* lies inside the boundary of atleast one of the circles

$$|\lambda - a_{kk}| = P_k \qquad (1)$$

Suppose above eq.(1) is for row-wise matrix, then similarly, using Gerschgorin theorem, we can write equation for column wise matrix. The intersection region gives the actual eigenvalues of the matrix *A*.

## III. BISECTION METHOD [ 3 ]

This is one of the simplest iterative methods and it is based on the property of intervals. To find a root using this method, let the function $f(x)$ be continuous between *a* and *b*. Suppose $f(a)$ is positive and $f(b), b = a + t, t > 0$ is negative. Then there is a root of $f(x) = 0$ lying between *a* and *b*.

## III. PROPOSED APPROACH FOR DETERMINING COMMON EIGENVALUES OF THE MATRICES

Suppose there are two matrices *A* and *B*. We need to determine common eigenvalues of above matrices. The various steps involved in determining common eigenvalues of the matrices as follows.

Step 1: Drawing Gerschgorin circles of matrices *A* and *B*.

Step 2: Determining intersection region of two matrices.

Step 3: Based on intersection region, determining bounds on the real axis in the s-plane.

Step 4: Using these bounds, under which all common eigenvalues are lying, we calculate common eigenvalues using Bisection method approach as discussed above.

### III. NUMERICAL EXAMPLE

Consider two matrices as

$$[A] = \begin{bmatrix} 3 & 1 & 4 \\ 0 & 2 & 6 \\ 0 & 0 & 5 \end{bmatrix}$$

and

$$[B] = \begin{bmatrix} 3 & -1 & 0 \\ -1 & 2 & -1 \\ 0 & -1 & 3 \end{bmatrix}$$

Now, we need to determine common eigenvalues using Gerschgorin theorem and Bisection method as follows.

**Proposed approach:**

*Step 1:* Drawing Gerschgorin circles for matrix $[A]$ and $[B]$ as shown in fig.1 and fig, 2.

Fig. 1.

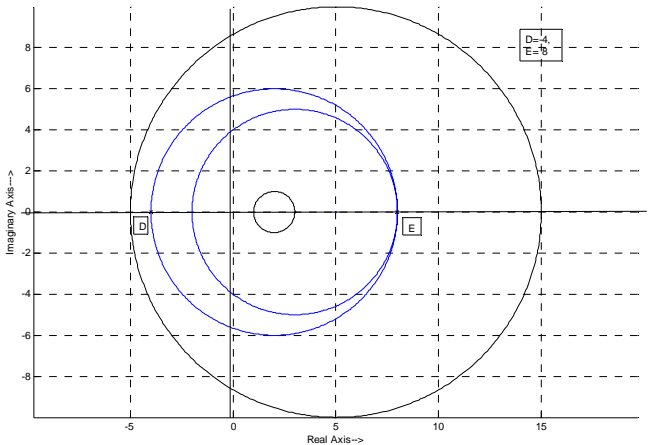

Fig 2.

The Gerschgorin bounds for matrix A is $D = -4, E = 8$.

The Gerschgorin bounds for matrix A is $D = 0, E = 4$.

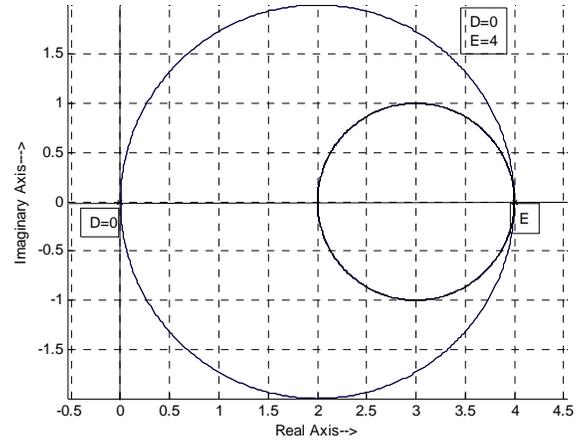

Step 2: Intersection of Gerschgorin circles of matrix *A* and *B* is shown in Fig. 3.

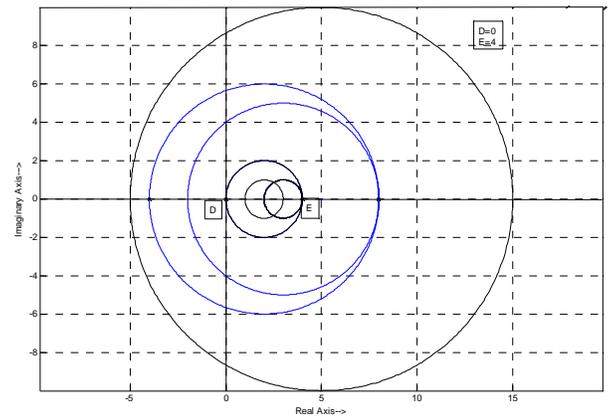

Fig. 3.

Step 3: The intersection region of bounds of matrix [A] and matrix [B] are $D = 0, E = 4$.

Step 4: Using Bisection method, common eigenvalues are determined by moving on real axis from $D = 0$ to $E = 4$ for both matrices $[A]$ and $[B]$. The results are shown in table 1 for matrix A and table 2 for matrix B.

| Sr. No. | $\lambda$ | $|\lambda I - A|$ | Remarks |
|---|---|---|---|
| 1 | 0 | -30 | |
| 2 | 0.1 | -26.9990 | |
| 3 | 0.2 | -24.192220 | |
| : | : | | |
| 21 | 1.9 | -1.3410 | Root=2 |
| 22 | 2 | 0 | |
| : | : | : | |
| 31 | 2.9 | 0.1890 | Sign change |
| 32 | 3 | -8.8818e-016 | Root $\approx 3$ |
| 42 | 4 | - 2.0790 | |

Table 2.

| Sr. No. | $\lambda$ | $|\lambda I - A|$ | Remarks |
|---|---|---|---|
| 1 | 0 | -12 | |
| 2 | 0.1 | -10.1790 | |
| 3 | 0.2 | -6.9930 | |
| : | : | : | |
| 10 | 0.9 | -0.6510 | Root=1 |
| 11 | 1 | 0 | |
| : | : | : | |
| 30 | 2.9 | 0.2090 | Sign change |
| 31 | 3 | -8.8818e-016 | Root $\approx 3$ |
| 41 | 4 | - 0.3410 | |

**Conventional approach:** using conventional method, bound for matrix $[A]$ is $D = -4, E = 8$. and matrix $[B]$ is $D = 0, E = 4$. Using bisection method, common eigenvalues are determined by moving on real axis from $D = -4$ to $E = 8$ for both matrices $[A]$ and moving on real axis from $D = 0$ to $E = 4$ for matrix $[B]$. The comparative results are shown in table 3.

Table 3.

| Sr. No | Methods | Computation time using Matlab |
|---|---|---|
| 1 | Conventional Method | 0.25 sec |
| 2 | Proposed approach | 0.016 sec |

Time required: 0.016 sec

Time required: 0.25 sec.

**Conclusions:**

It is observed that in the proposed method the algorithm takes much time then the conventional method to compute the common eigenvalues between the two matrices. The common eigenvalues are used in the identification of two images during the image recognition processes. It is also useful in Control theory and Computer Engineering Applications.

ACKNOWLEDGEMENT

We are thankful to Dr. Yogesh V. Hote for suggestion in writing this paper.